\documentclass[12pt]{article}
\usepackage{theorem,amsfonts,amssymb}

\newtheorem{theorem}{Theorem}[section]
\newtheorem{proposition}{Proposition}[section]
\newtheorem{lemma}{Lemma}[section]
\newtheorem{corollary}{Corollary}[section]
\theorembodyfont{\upshape}
\newtheorem{definition}{Definition}[section]
\newtheorem{remark}{Remark}[section]
\newtheorem{example}{Example}[section]
\newtheorem{proof}{Proof}

\newcommand{\bt}{\begin{theorem}}
\newcommand{\et}{\end{theorem}}
\newcommand{\bl}{\begin{lemma}}
\newcommand{\el}{\end{lemma}}
\newcommand{\bp}{\begin{proposition}}
\newcommand{\ep}{\end{proposition}}
\newcommand{\bd}{\begin{definition}}
\newcommand{\ed}{\end{definition}}
\newcommand{\br}{\begin{remark}}
\newcommand{\er}{\end{remark}}
\newcommand{\bex}{\begin{example}}
\newcommand{\eex}{\end{example}}
\newcommand{\bc}{\begin{corollary}}
\newcommand{\ec}{\end{corollary}}
\newcommand{\bo}{\begin{proof}}
\newcommand{\eo}{\end{proof}}
\newcommand{\be}{\begin{enumerate}}
\newcommand{\ee}{\end{enumerate}}

\title{Liouville property on $G$-spaces}
\author{C. R. E. Raja}
\date{ }

\begin{document}
\maketitle

\let\epsi=\epsilon
\let\vepsi=\varepsilon
\let\lam=\lambda
\let\Lam=\Lambda
\let\ap=\alpha
\let\ol=\overline
\let\vp=\varphi
\let\ra=\rightarrow
\let\Ra=\Rightarrow
\let\da=\downarrow
\let\Llra=\Longleftrightarrow
\let\Lla=\Longleftarrow
\let\lra=\longrightarrow
\let\Lra=\Longrightarrow
\let\ba=\beta
\let\ga=\gamma
\let\Ga=\Gamma
\let\un=\upsilon
\let\ct=\cdot
\let\te=\tilde

\newcommand{\cI}{{\cal I}}
\newcommand{\cSL}{{\cal SL}}
\newcommand{\cG}{{\cal G}}
\newcommand{\cA}{{\cal A}}
\newcommand{\mP}{{\mathbb P}}
\newcommand{\N}{{\mathbb N}}
\newcommand{\C}{{\mathbb C}}
\newcommand{\cX}{{\cal X}}
\newcommand{\cY}{{\cal Y}}
\newcommand{\R}{{\mathbb R}}

\begin{abstract}
Let $G$ be a locally compact group and $E$ be a $G$-space.  An
irreducible probability measure $\mu$ on $G$ is said to have
Liouville property on $E$ if $G$-invariant functions on $E$ are the
only continuous bounded functions on $E$ that satisfy the mean value
property with respect to $\mu$.  We first prove that the random walk
induced by $\mu$ on $E$ is transient outside a closed set and on the
closed set $\mu$ has Liouville.  We mainly consider actions on
vector spaces and projective spaces.  We show that measures on
$GL(V)$ that are supported inside a ball of radius less than $a<1$
have Liouville property on $V$.  We also prove that measures on
$GL(\R ^2)$ have Liouville property on the projective line.  We next
exhibit subgroups of $GL(V)$ so that irreducible measures on such
subgroups have Liouville on the projective space $\mP (V)$ of $V$.
We also prove irreducible measures on $SL(V)$ have Liouville
property on $\mP (\cSL (V))$ where $\cSL (V)$ is the Lie algebra of
$SL(V)$.
\end{abstract}

\begin{section}{Introduction}

Let $G$ be a locally compact group and $\mu$ be a regular Borel
probability measure on $G$.  We consider a locally compact space $E$
on which the group $G$ acts by homeomorphisms. A bounded measurable
function on $E$ that satisfies the mean value property with respect
to $\mu$ is called a $\mu$-harmonic function and $H_\mu (E)$ is the
space of all bounded $\mu$-harmonic functions on $E$. In case of
left or right action of $G$ on itself, the space of harmonic
functions was introduced by Fursternberg \cite{Fu1} and studied by
others. An earlier work of Blackwell, Choquet and Deny on abelian
groups showed that constant functions are the only continuous
bounded $\mu$-harmonic functions on abelian groups $G$ for the left
or right action of $G$ on itself - such a result is known as
Choquet-Deny theorem or Liouville property (cf. \cite{JaR} for
recent developments in Choquet-Deny results on groups).  Here, we
consider Liouville property for group actions. In this situation we
say that $\mu$ has Liouville property on a $G$-space $E$ if
$G_\mu$-invariant functions are the only continuous functions in
$H_\mu (E)$ where $G_\mu$ is the closed subgroup generated by the
support of $\mu$.

By considering the adjoint $\check \mu$ of $\mu$, it may be seen
that $\check \mu$ has Liouville on $G$ implies $\mu$ has Liouville
on any $G$-space $E$: recall $\check \mu$ is defined by $\int f(g)
d\check \mu (g)= \int f(g^{-1}) d\mu (g)$.  However we provide
examples in Section 7 to show that there are measures on $GL(V)$
that have Liouville property on $V$ but neither the measure nor its
adjoint has Liouville property on $GL(V)$.

It may be easily observed that there is a $S_\mu$-invariant closed
(possibly empty) subset $L _\mu$ of $E$ such that any continuous
function in $H_\mu (E)$ is $S_\mu$-invariant when restricted to
$L_\mu$: $S_\mu$ is the closed semigroup generated by the support of
$\mu$.

Recently, R. Feres and E. Ronshausen \cite{FeR} studied Liouville
property for group actions and it was shown that if $\Ga$ is a
countable group acting on the circle $S^1$ or $[0,1]$, then any
irreducible (that is, $S_\mu =G$) symmetric measure on $\Ga$ has
Liouville property on $S^1$ or $[0,1]$.

It is also shown in \cite{FeR} that the random walk generated by
irreducible $\mu$ is transient on $E\setminus L_\mu$ for actions of
countable groups on compact spaces using boundaries of countable
groups, we extend this result to all type of actions using results
about harmonic functions on groups (see Proposition \ref{hd}).

We mainly consider actions on vector spaces and on the corresponding
projective spaces. For a finite-dimensional vector space $V$ over
$\R$, and a linear transformation $\ap $ of $V$, $||\ap ||$ denotes
the operator norm of $\ap$.  We first observe the following result
which serves many interesting examples.

\bt\label{gc} Let $V$ be a finite-dimensional vector space over $\R$
and $G$ be a subgroup of $GL(V)$ such that $\{ g \in G \mid ||g|| <1
\}$ is a non-empty open set.  Suppose $\mu$ is an adapted
probability measure on $G$ such that $\mu ( \{ g\in G \mid ||g||
\leq  a \}) =1 $ for some $0<a <1$.  Then $\mu$ has Liouville
property on $V$. \et

For a $G$-space $E$, we say that $G$ has Liouville property on $E$
if any irreducible probability measure on $G$ has Liouville on $G$.

We next consider projective linear actions on projective spaces. For
a vector space $V$, $\mP (V)$ denotes the corresponding projective
space.  We first prove Liouville property for actions on the
projective line $\mP ^1= \mP (\R ^2)$: it may be recalled that $S^1$
and the projective line $\mP ^1= \mP (\R ^2)$ are homeomorphic and
\cite{FeR} proved Liouville property for (symmetric measures on)
countable group actions on $S^1$.

\bt\label{dim2} Let $G$ be a closed subgroup of $GL( \R ^2)$.  Then
$G$ has Liouville property on $\mP ^1$. \et

We next look at locally compact subgroups of $GL(V)$.

\bt\label{ur} Let $G$ be locally compact subgroup of $GL(V)$.
Suppose $G$ has an unipotent subgroup $U$ such that $U$ has only one
linearly independent invariant vector in $V$.  Then $G$ has
Liouville property on $\mP (V)$. \et

\bex  The following Lie subgroups of $GL(V)$ have unipotent
subgroups that have only one linearly independent invariant vector
in $V$.

\be

\item $G$ is the group of invertible upper triangular matrices and $U$ may be
taken as the group of upper triangular matrices with one on the
diagonal.

\item the $(2n+1)$-dimensional Heisenberg group and any of its
extensions.  For instance, $SL_2(\R ) \ltimes H_1$ with $$H_1 = \{
\pmatrix { 1&a &x\cr 0&1&b \cr 0&0&1 } \mid a,x,b \in \R \}$$ as the
unipotent group.

\ee \eex

As an illustration and application of the method involved in the
proof of Theorem \ref{ur}, we consider the conjugate action of
$SL(V)$ on its Lie algebra $\cSL (V)$. It may be noted that the Lie
algebra $\cSL (V)$ of $SL(V)$ is the space of all trace zero
matrices in $GL(V)$ and the adjoint action is given by ${\rm Ad}~(g)
(x) = gxg^{-1}$ for all $g\in SL(V)$ and $x\in \cSL (V)$.

\bt\label{ca} Let $V$ be a finite-dimensional vector space over
$\C$.  For the adjoint action of $SL(V)$ on the Lie algebra $\cSL
(V)$, $SL(V)$ has Liouville property on $\mP (\cSL (V))$. \et

As a by-product we obtain that actions such as distal, proximal and
minimal have Liouville property (see Corollary \ref{prx}). We make a
few miscellaneous remarks on amenability of the group $G$ and
Liouville property of $G$-spaces and regarding equicontinuous Markov
operator associated to the random walk on $E$ generated by $\mu$.

\end{section}

\begin{section}{Preliminaries}

Let $G$ be a locally compact group and $\mu$ be a (regular Borel)
probability measure on $G$.

Let $G_\mu$ (resp. $S_\mu$) denote the closed subgroup (resp.
semigroup) generated by the support of $\mu$. The measure $\mu$ is
called adapted (resp. irreducible) if $G_\mu =G$ (resp. $S_\mu =G$).

A $G$-space $E$ is a locally compact space with a continuous
$G$-action.

A bounded function $f$ on $E$ is called $\mu$-harmonic if $Pf(x)=
\int f(gx)d\mu (g)=f(x)$ for all $x\in E$.

Let $H_\mu (E)$ denote the space of bounded $\mu$-harmonic functions
on $E$ and $C_b(E)$ denote the space of bounded continuous functions
on $E$.

We say that the measure $\mu$ has Liouville property (or
Choquet-Deny) on $E$ if any continuous bounded $\mu$-harmonic
function on $E$ is $G_\mu$-invariant.  We say that (the action of)
$G$ has Liouville property on $E$ if every irreducible probability
measure on $G$ has Liouville property on $E$.

For right-translation action of $G$ on itself, we simply say that
$\mu$ has Choquet-Deny if any continuous bounded right
$\mu$-harmonic function is constant on the cosets of $G_\mu$.  It is
easy to verify that if $\mu$ has Choquet-Deny, then $\check \mu$ has
Liouville property on any $G$-space: $\check \mu $ is given by
$\check \mu (E) = \mu (E^{-1})$ where $E ^{-1}= \{ x^{-1} \mid x\in
E \}$.

\begin{subsection}{Recurrence and Transience}

Let $(X_n)$ (resp. $(Z_n)$) be the canonical left (resp. right)
random walk on $G$ defined by $\mu$. Then for $x\in E$, $(X_nx)$
defines a random walk on $E$ starting at $x$.

For a given $x\in E$, we consider the following recurrence property
$R_x$ and transition property $T_x$ introduced in \cite{GR} for
$(X_nx)$:

$R_x$: There exists a compact set $K\subset E$ such that a.e.
$X_nx\in K$ infinitely often.

$T_x$: For any compact set $K \subset E$, a.e. $\omega$, there
exists $n(w)\in \N$ with $X_n(\omega )x \not \in K$ for $n \geq
n(\omega )$.

Several situations for the validity of $R_x$ or $T_x$ are given in
\cite{GR}.

For the action of $G$ on $E$ and a probability measure $\mu$ on $G$,
we have a closed $S_\mu$-invariant subset $L$ of $E$ such that any
continuous $f\in H_\mu (E)$, $f$ is $S_\mu$-invariant on $L$.  It is
shown in \cite{FeR} that the random walk is transient on $E\setminus
L$ using Poisson boundary, here we first extend this result to
general actions using a well-known result about harmonic functions.

\bp\label{hd} There is a $S_\mu$-invariant closed set $L$ such that
$f$ is $S_\mu$-invariant on $L$ for any continuous $f\in H_\mu (E)$
and the limit points of almost all paths $(X_nx)$ are in $L$ for any
$x\in E$.  The action of $G$ on $E\setminus L$ satisfies $T_x$ for
each $x\in E\setminus L$.  In other words, if $R_x$ is valid for all
$x\in E$, then any bounded continuous $\mu$-harmonic function on $E$
is $G_\mu$-invariant. \ep

The following lemma easily follows from a similar well-known result
about $\mu$-harmonic functions on groups, we include the proof for
clarity and completeness.

\bl\label{mc} For any any $f\in H_\mu (E)$ and $x\in E$, we have for
$\mu$-almost every $g\in G$, $\lim f(gX_nx) - f(X_nx) =0$ a.e.
$\omega$. \el

\bo Let $F(g) = f(g^{-1}x)$.   Then $F$ is a (right) $\check
\mu$-harmonic function on $G$.  Let $(X_n)$ be the left random walk
induced by $\mu$.  Then $(\tilde X_n=X_n^{-1})$ is the right random
walk induced by $\check \mu$.  It could be easily verified that
$F(\tilde X_n)$ is a martingale, that is $E(F(\tilde X_{n+1})|\cA
_n) =F(\tilde X_n)$ where $\cA _n$ is the $\sigma$-algebra generated
by $F(\tilde X_1), \cdots , F(\tilde X_n)$. Since $F(\tilde X_n)$ is
bounded,
$$\begin{array}{ll} E((F(\tilde X_{n+1})-F(\tilde X_n))^2) & \\
= E(E((F(\tilde X_{n+1})-F(\tilde X_n))^2|\cA _n)) & \\ =
E(E(F(\tilde X_{n+1})^2| \cA _n))- 2E(F(\tilde X_n)E(F(\tilde
X_{n+1})|\cA _n)) +E(F(\tilde X_n)^2) & \\ = E(F(\tilde X_{n+1})^2)-
E(F(\tilde X_n)^2)
\end{array}$$ for all $n$.  This implies that $\sum _{k=0}^n
E((F(\tilde X_{k+1})-F(\tilde X_k))^2) = E(F(\tilde X_{n+1}^2)) -
E(F(\tilde X_0^2))$. Since $E(F(\tilde X_n)^2)$ is a uniformly
bounded, by martingale convergence theorem $F(\tilde X_n)$ converges
almost surely and in $L^2$, hence $E(F(\tilde X_{n+1}^2))$
converges.  This implies that $\sum E((F(\tilde X_{n+1})-F(\tilde
X_n))^2) <\infty$ and hence $\sum E((F(\tilde X_ng^{-1})-F(\tilde
X_n))^2) <\infty$ for $\mu$-almost every $g$. This implies that
$\sum (F(\tilde X_ng^{-1})-F(\tilde X_n))^2 <\infty$ a.e. $\omega$.
Thus, $F(\tilde X_ng^{-1})-F(\tilde X_n)\to 0$ a.e. $\omega$. \eo

\bo $\!\!\!\!\!$ {\bf of Proposition \ref{hd}} Let $L$ be the
largest subset of $E$ consisting of all $x\in E$ such that $f(x) =
f(gx)$ for all $g\in S_\mu$ and all continuous $f\in H_\mu (E)$.
Then $L$ is a closed $S_\mu$-invariant subset of $E$.

If $X_n$ is the left random walk on $G$ defined by $\mu$ and $f\in
H_\mu (E)$, then by Lemma \ref{mc}, for a fixed $x\in E$,
$$ \lim f(gX_nx) - f(X_nx) =0 \eqno (1)$$ almost all $\omega$ and
$\mu$-a.e. $g$.  If $y$ is a limit point of $(X_nx)$ for a path
$(X_nx)$ which satisfies (1), then $f(gy)=f(y)$ for all $g\in S_\mu$
and all continuous $f\in H_\mu (E)$.  This implies that the limit
points of almost all paths are in $L$.

Let $x\in E\setminus L$ and $K$ be any compact subset of $E\setminus
L$.  If $X_n(\omega )x\in K$ infinitely often, then $X_n(\omega)x$
has a limit point $y\in K$.  By the first part such $\omega $ form a
null set.  Hence there is an almost finite random integer $N_x$ such
that $X_n(\omega )x\not \in K$ for all $n \geq N(\omega )$. \eo

\end{subsection}

\begin{subsection}{A sufficient condition}

We now provide a sufficient condition useful to prove the main
results.

Following result is a well-known useful basic result in harmonic
functions, as it is a standard result which is often needed we state
it without proof.

\bl\label{ls} Let $f$ be a continuous bounded $\mu$-harmonic
function on $E$.  Then the sets $\{ x\in E \mid f(x) = \inf _{a\in
E} f(a) \}$ and $\{ x\in E \mid f(x) = \sup _{a\in E} f(a) \}$ are
$S_\mu$-invariant closed sets (possibly empty). \el

Liouville property for $\mu$ rests in showing that the two
$S_\mu$-invariant sets in the above Lemma \ref{ls} intersect
(instead of $E$, one considers $\ol{G_\mu x}$).  This motivates us
to provide the following sufficient condition for Liouville
property-it is a transparent version of Corollary 2.2 of \cite{FeR}.

\bp\label{um}  If $f$ is a continuous $\mu$-harmonic function on $E$
such that $f$ has a minimum and maximum on $E$ and $S_\mu$-invariant
sets overlap, then $f$ is constant.  In particular, if for any $x\in
E$, $\ol{Gx}$ is compact and $S_\mu$-invariant subsets of $\ol {Gx}$
overlap, then $G$ has Liouville property on $E$. \ep

\br The condition that $S_\mu$-invariant subsets of $\ol {Gx}$
overlap is not a necessary condition.  Example \ref{ex1} shows that
there are actions having Liouville property but have orbits
violating this condition.  In section 5, we provide a class of
measures/actions for which the above sufficient condition is
necessary. \er

\bo Let $f$ be a continuous $\mu$-harmonic function having maximum
and minimum on $E$.  Let $E_s = \{ x \in E \mid f(x) = \sup _{y\in
E} f(y) \}$ and $E_i = \{ x \in E  \mid f(x) = \inf _{y\in E} f(y)
\}$. Then $E_s$ and $E_i$ are nonempty closed subsets in $E$.  By
Lemma \ref{ls}, both $E_s$ and $E_i$ are $S_\mu$-invariant closed
subsets of $E$. If $S_\mu$-invariant sets pairwise overlap, then
there is a $x\in E_s\cap  E_i$ and hence $f$ is constant on $E$.
Apply the first result to each orbit closure to obtain the second
result.  \eo

As a consequence we get the Liouville property for proximal actions,
minimal actions and distal actions on compact spaces.  Recall that a
semigroup $S$ acting a compact space $E$ is called proximal (resp.
distal) if for any two distinct points $x,y \in E$, the closure of
$\{ (gx, gy) \mid g\in S \}$ meets (resp. does not meet) the
diagonal in $E\times E$.

\bc\label{prx} Let $G$ be a locally compact group acting on a
compact space $E$ and $\mu$ be a probability measure on $G$.

\be

\item[1] If the action of $S_\mu$ on $E$ is orbitwise proximal, that is proximal
on closure of any orbit, then $\mu$ has Liouville property on $E$.

\item[2] If the orbit closures $\overline {S_\mu x}$ are minimal, then $\mu$  has Liouville
property on $E$.

\item[3] If the action of $S_\mu$ on $E$ is distal, then $\mu$  has Liouville
property on $E$. \ee \ec

\bo  We assume that $E= \ol {S_\mu a}$ for some $a\in E$.

If the action is orbitwise proximal, then $S_\mu$ is proximal on
$E$.  Since $E$ is compact, $\ol{S_\mu x} \cap \ol{S_\mu y}\not =
\emptyset$ for any $x, y \in E$.  This implies that any two
$S_\mu$-invariant sets overlap. Now the result follows from
Proposition \ref{um}.

If the orbit closures are minimal, then the result easily follows
from Proposition \ref{um}.

If the action is distal, then the orbit closures are minimal, hence
the result follows from 2. \eo

\end{subsection}
\end{section}

\begin{section}{Actions on vector spaces}

We now look at actions on vector spaces.  Given a sequence $(\mu
_n)$ of probability measures on a locally compact space $E$, we say
that $\mu _n \to \mu $ in the weak* topology for a probability
measure $\mu$ on $E$ if $\mu _n(f) \to \mu (f)$ for any continuous
bounded function $f$ on $E$.  It may be noted that $\mu _n \to \mu$
in the weak* topology if and only if $\mu _n(f) \to \mu (f)$ for any
continuous function $f$ with compact support on $E$ (cf. 1.1.9 of
\cite{He}).

\bo $\!\!\!\!\!$ {\bf of Theorem \ref{gc}} Since $\mu ( \{ g\in G
\mid ||g ||\leq  a \} )=1$ for some $0< a <1$, $\mu ^n ( \{ g\in G
\mid ||g ||\leq  a^n \}) =1$ for $n \geq 1$.  Let $v\in V$.  Then
$||gv|| \leq a^n ||v|| $ for all $g$ in the support of $\mu ^n$. Let
$\epsi
>0$ and $\psi$ be any continuous function with compact support on
$V$.  Then using uniform continuity of $\psi$, we get that $|\psi
(gv)-\psi (0)| <\epsi$ for all $g$ in the support of $\mu ^n$, for
large $n$.  This implies that $\mu ^n
*\delta _v (\psi ) = \int \psi (gv) d\mu ^n (g)\to \psi (0) = \delta _0
(\psi )$, hence $\mu ^n *\delta _v \to \delta _0$ in the weak*
topology.

Let $f$ be a continuous bounded $\mu$-harmonic function on $V$. Then
$f(v) = \mu ^n *\delta _v(f) \to f(0)$ for any $v\in V$, hence $f$
is constant.  Thus, $\mu$ has Liouville property on $V$. \eo

\end{section}

\begin{section}{Actions on projective spaces}

We now consider actions on projective spaces.  We first prove a
useful result on unipotent actions on projective spaces.  A split
solvable algebraic group is a solvable algebraic group whose maximal
torus is a split torus: e.g. unipotent algebraic groups, $(\R ^*)^n
\ltimes \R^n$ where $\R ^* = \R \setminus \{ 0 \}$, more generally,
group of all upper triangular matrices.

\bl\label{ss} Let $V$ be a finite-dimensional vector space over $\R$
and $G$ be a closed subgroup of $GL(V)$.  If $G$ is amenable, then
action of $G$ on any $G$-minimal subset of $\mP (V)$ factors through
a compact group.  If the algebraic closure of $G$ is a split
solvable algebraic group, then any $G$-minimal subset of $\mP (V)$
is singleton consisting of a $G$-fixed point. \el

\bo  Let $E$ be a $G$-minimal subset of $\mP (V)$.  If $G$ is
amenable, then $E$ supports a $G$-invariant probability measure
$\rho$. Since $E$ is $G$-minimal, support of $\rho$ is $E$.  Let
$I_\rho= \{ g\in GL(V) \mid g(x) =x {\rm ~~for ~~all ~~} x\in E \}$
and $\cI _\rho = \{ g\in GL(V) \mid g(\rho ) = \rho \}$. Then
$G\subset \cI _\rho$ and $\cI _\rho /I_\rho $ is a compact algebraic
group (cf. \cite{Da} and \cite{Fu2}). If the algebraic closure of
$G$ is a split solvable algebraic group, $G\subset I_\rho$.\eo

We next consider actions on the projective line.

\bo $\!\!\!\!\!$ {\bf of Theorem \ref{dim2}} Let $\phi \colon \R ^2
\setminus \{ 0\} \to \mP ^1$ be the canonical projection. Suppose
$\R^2$ is not $G$-irreducible. Then there is a nonzero vector $v\in
\R ^2$ such that for each $g\in G$, $gv=t_gv$ for some $t_g\in \R$.
If all $g\in G$ is diagonalizable over $\R$, then $G$ is abelian,
hence $G$ has Liouville on $\mP ^1$. So, assume that there is a
$g_0\in G$ that is not diagonalizable over $\R$. Thus, $v$ is the
only eigenvector of $g _0\in G$.  Let $G_0$ be the group generated
by $g_0$. Then $\phi (v)$ is the only $G_0$-invariant vector and the
algebraic closure of $G_0$ is a split solvable algebraic group.  By
Lemma \ref{ss}, any $G_0$-minimal subset of $\mP ^1$ consists of a
$G_0$-fixed point. Since $\phi (v)$ is the only $G_0$-fixed point,
we get that any $G$-invariant closed subset of $\mP ^1$ contains
$\phi (v)$.  Now the result follows from Proposition \ref{um}.

We now assume that $\R^2$ is $G$-irreducible.  Suppose $G$ is not a
relatively compact subset of projective linear transformation of
$\mP ^1$, that is, $PGL (\R ^2)$.  Then there is a sequence $(g_n)$
in $G$ such that $(g_n)$ has no convergent subsequence in $PGL(\R
^2)$. Passing to a subsequence, we may assume that ${g_n \over
||g_n||} \to h$ where $h$ is a linear transformation on $\R ^2$.
Since $||h ||=1$, $h\not = 0$.  Since $(g_n)$ has no convergent
subsequence in $PGL(\R ^2)$, $h$ is a rank-one transformation.  Let
$w$ be a nonzero vector in the image of $h$. Then $g_n (\phi (x))
\to \phi (w)$ for all $x$ not in the kernel of $h$.  If $x$ is a
nonzero vector in the kernel of $h$, then since $\R^2$ is
$G$-irreducible, there is a $g\in G$ such that $gx$ is not in the
kernel of $h$.  Now ${g_ng\over ||g_ng||}\to {hg \over ||hg||}$.
Since $h(gx)\not =0$ we get that $g_ng(\phi (x)) \to \phi (hg(x))=
\phi (w)$.  Thus, $\phi (w) \in \ol {G\phi (x)}$ for any non-zero
$x\in \R^2$.  Now the result follows from Proposition \ref{um}. \eo

We next consider actions of locally compact subgroups of $GL(V)$.

\bo $\!\!\!\!\!$ {\bf of Theorem \ref{ur}} Let $U$ be an unipotent
subgroup of $G$ that has only one linearly independent invariant
vector $v \in V$.  Suppose $x \in \mP(V)$ is $U$-invariant.  Let
$\phi \colon V\setminus \{ 0 \} \to \mP(V)$ be the canonical
projection and $w\in V\setminus \{ 0 \}$ be such that $\phi (w) =x$.
Then for each $g\in U$ there is a $t_g \in \R$ such that $g(w)
=t_gw$. Since $g\in U$ are unipotent, $g(w)=w$.  Since $U$ has only
one linearly independent invariant vector $v \in V$, $w= tv$ for
some $t\in \R$, that is, $\phi (v) = \phi (w)$.

Let $E$ be a $G$-invariant closed subset of $\mP (V)$.  Then $E$ is
$U$-invariant and contains a $U$-minimal subset $M$.  By Lemma
\ref{ss}, $M$ consists of a fixed point of $U$, hence $M= \{ \phi
(v) \}$. Thus, any $G$-invariant closed subset contains $\phi (v)$.
This implies by Proposition \ref{um}, that any continuous bounded
$\mu$-harmonic function is constant. \eo

We now look at the conjugate action of $SL(V)$.

\bo $\!\!\!\!\!$ {\bf of Theorem \ref{ca}} Let $e_{1,n}$ be the
matrix whose $(i,j)$-th entry is nonzero only for $i=1$ and $j=n$
and $\phi \colon \cSL (V) \setminus \{ 0\} \to \mP(\cSL (V))$ be the
canonical projection.  We now prove that every $SL(V)$-invariant
subset of $\mP(\cSL (V))$ contains $\phi (e_{1,n})$.

Choose subspaces $\{0 \}=W_0\subset W_1\subset \cdots \subset W_d=V$
such that ${\rm dim}(W_i)=i$.  Take $$U= \cap _{1\leq k <d } \{ g\in
SL(V) \mid g(W_k)= W_k ~~{\rm and } ~~ g(v)-v\in W_k ~~{\rm for ~~
all }~~v\in W_{k+1} \}.$$  Then $U$ is the unipotent algebraic group
of all upper triangular unipotent matrices in $SL(V)$ (with respect
to a basis from the subspace $W_k$). It is easy to notice that
$W_{i+1}/W_i$ is the subspace of all $U$-invariant vectors in
$W_d/W_i$ for all $0\leq i <d$ and can easily be seen that the
center of $U$ is $\{ \exp (te_{1,n}) \mid t\in \R \}$.

Let $v \in \cSL(V)$ be such that $v$ is $U$-invariant, that is,
$gv=vg$ for all $g\in U$.  Let $v= v_s+v_n$ be the Jordan-Chevalley
decomposition of $v$ into a semisimple element $v_s$ and a nilpotent
element $v_n$ such that $v_sv_n= v_nv_s$ and there are polynomials
$P$ and $Q$ with $v_s = P(v)$ and $v_n = Q(v)$ (cf. Proposition 4.2
of \cite{Hu}). Then since $gv=v$ for all $g\in U$, $gv_s=v_sg$ and
$gv_n =v_ng$ for all $g\in U$.   Since $v_n u=uv_n$ for all $u\in
U$, any $U$-invariant vector is also $v_n$-invariant, hence $W_1$ is
invariant under $v_n$. Since $v_n$ is nilpotent and ${\rm dim}
(W_1)=1$, $v_n(W_1)= \{0 \}=W_0$. We now claim by induction that
$v_n(W_k)\subset W_{k-1}$ for all $k\geq 1$. Suppose
$v_n(W_i)\subset W_{i-1}$ for some $1\leq i <d$.  Since $v_n u=uv_n$
for all $u \in U$, by considering the $U$-invariant one-dimensional
subspace $W_{i+1}/W_i$ of $W_d/W_i$, we conclude that
$v_n(W_{i+1}/W_i)\subset W_{i+1}/W_i$.  Since $v_n$ is nilpotent and
${\rm dim}(W_{i+1}/W_i )=1$, we get that $v_n(W_{i+1})\subset W_i$.
This implies that $\exp {(v_n)}\in U$. Since $v_n$ commutes with
every element of $U$, $\exp {(v_n)}$ is in the center of $U$. Thus,
$v_n =te_{1,n}$ for some $t\in \C$. Since $uv_s=v_su$ for all $u \in
U$, it can easily be seen that $v_s$ has only one eigenvalue and
hence $v_s =tI$ for some $t\in \C$.  Since $v$ has trace zero,
$v_s=0$. Thus, $v= v_n=te_{1,n}$. This shows that $U$ has only one
linearly independent invariant vector $e_{1,n}$ in $\cSL(V)$.  Now
the result may be proved as in Theorem \ref{ur}.\eo

\end{section}

\begin{section}{Miscellaneous remarks}

\begin{subsection}{Amenability and Liouville property}

A well-known conjecture of Furstenberg \cite{Fu1} states that a
locally compact $\sigma$-compact group $G$ is amenable if and only
if $G$ admits adapted probability measures having Liouville
property.  Several proofs of this conjectures are available (cf.
\cite{KaV}, \cite{Ro} and \cite{Wi}).  In this case one could prove
a similar result for Liouville property of group actions using the
above conjecture.

\bp\label{ca} Let $G$ be a locally compact $\sigma$-compact group
and $\mu$ be an adapted probability measure on $G$.  Then the
following are equivalent:

\be \item [(1)] $\check \mu$ has Choquet-Deny;

\item [(2)] $\mu$ has Liouville property on all
$G$-spaces;

\item [(3)] $\mu$ has Liouville property on all compact
$G$-spaces;

\item [(4)] $\mu$ has Liouville property on all compact affine
$G$-spaces. \ee

In particular, $G$ is amenable if and only if $G$ has an adapted
probability measure $\mu$ such that $\mu$ has Liouville property on
all compact $G$-spaces. \ep

\bo It is sufficient to prove that (4) $\Ra$ (1).  Suppose $\mu$ has
Liouville property on all compact affine $G$-spaces.  Let $C_b(G)$
be the Banach space of all continuous bounded functions on $G$. Let
$E$ be the unit ball in the dual $C_b(G)^*$ of $C_b(G)$: it may be
noted that the dual of $C_b(G)$ is the space of all regular bounded
finitely additive measures on $G$.  Then $E$ is compact in the weak*
topology.  $G$ acts on $C_b(G)$ by $gf(x) = f(xg)$ for all $x,g \in
G$ and $f\in C_b(G)$ and the action is by isometries.  By duality
each $g\in G$ defines a continuous map $g^*$ on $E$.  Now, $g
\mapsto {g^{-1}}^*$ defines an action of $G$ on $E$.

Let $f$ be a continuous bounded $\check \mu$-harmonic function on
$G$, that is, $$ \int f(xg^{-1}) d\mu (g) = f(x), ~~ x\in G .$$
Define the function $\phi$ on $E$ by $\phi (\sigma ) = <\sigma, f>$
for $\sigma \in E$.  Then $\phi$ is a continuous bounded function on
$E$.  Now, $\int \phi ({g^{-1}}^*\sigma ) d\mu (g)= \int
<{g^{-1}}^*\sigma, f> d\mu (g) = \int <\sigma, g^{-1}f> d\mu (g) =
<\sigma , \int g^{-1}f d\mu(g)> =<\sigma,  f>= \phi (\sigma )$.
Thus, $\phi$ is a $\mu$-harmonic function on $E$.  By assumption
$\phi$ is constant on $G$-orbits, that is, $<{x^{-1}}^*\sigma, f>=
<\sigma , f>$ for all $x\in G$ and $\sigma \in E$.  For any $x\in
G$, since $\delta _x \in E$ we have $f(x) = \phi (\delta _x) =
\phi(x^{-1}\delta _e) = \phi (\delta _e)= f(e)$. Thus, $f$ is
constant.\eo

\end{subsection}

\begin{subsection}{Equicontinuous Markov operator}

If $E$ is a $G$-space and $\mu$ is a probability measure on $G$,
then we define a Markov operator $P$ on $C_b(E)$ by $$Pf(x) = \int
f(gx) d\mu (g) , ~~~ f\in C_b(E), ~~x\in E .$$  We say that the
Markov operator $P$ defined by $\mu$ on $E$ is equicontinuous if for
every $f\in C_b(E)$, there exists a subsequence of integers $(k_n)$
and $F\in C_b(E)$ such that $\{ {1\over k_n} \sum _{i=1} ^{k_n} P^i
f \} \to F$ pointwise.

\br  \be

\item If the closed semigroup (or equivalently the closed subgroup)
generated by the support of $\mu$ is compact, then $P$ is
equicontinuous on any $G$-space $E$ which may be seen as follows.
Since the closed semigroup generated by the support of $\mu$ is
compact, ${1\over n}\sum _{k=1} ^n \mu ^k $ converges to a
probability measure $\lam$ in the weak*-topology, that is, ${1\over
n}\sum _{k=1} ^n \mu ^k (f) \to \lam (f)$ for all continuous bounded
functions $f$ on $G$ (cf. \cite {Ro2}).  Thus, if $E$ is a
$G$-space, then for $x\in E$, ${1\over n} \sum _{k=1} ^{n} \mu ^k
*\delta _x \to \lam *\delta _x$ in the weak* topology on $E$,
hence ${1\over n} \sum _{k=1} ^{n} P^k f (x) \to \lam *\delta _x(f)$
for all continuous bounded function $f$ on $E$.  Thus, $P$ is
equicontinuous on $E$.

\item If $V$ is a finite-dimensional vector space and $\mu$ is a
measure on $GL(V)$ such that $S _\mu$ is strongly irreducible on
$V$, that is, no finite union of proper subspaces is invariant under
$S _\mu$, then by Proposition 3.1 of \cite{BQ} we get that $P$
defined by $\mu$ on $\mP (V)$ is equicontinuous.

\ee\er

We now characterize equicontinuous $P$ that has Liouville property.

\bt\label{ec}  Let $E$ be a compact $G$-space and $\mu$ be an
adapted probability measure on $G$ such that the corresponding
Markov operator $P$ on $E$ is equicontinuous. Then the following are
equivalent:

\be

\item[(1)] $\lam$ has Liouville property on $E$ for all adapted probability
measures $\lam$ on $G$;

\item[(2)] $\mu$ has Liouville property on $E$;

\item[(3)] for any $x\in E$, $\ol{S_\mu x}$ has a unique minimal subset;

\item[(4)] for any $x\in E$, $S_\mu$-invariant subsets of $\ol{Gx}$
overlap.\ee \et

\bo (1) $\Ra $ (2) is evident and (4) $\Ra$ (1) follows from
Proposition \ref{um}.  It only remains to show that (2) $ \Ra $ (3)
$\Ra $ (4).

If there is a $x \in E$ such that $\ol {S_\mu x}$ contains disjoint
nonempty closed $S_\mu$-minimal sets $E_1$ and $E_2$. Let $f$ be a
continuous function on $E$ such that $f(E_1) = 1$ and $f(E_2 ) =2$.
Since $P$ is equicontinuous, there exists $(k_n)$ such that ${1\over
k_n}\sum _{i=1}^{k_n} P^i f $ converges to a continuous function $F$
on $E$.  It can easily be verified that $PF=F$ but $F(E_1)=1$ and
$F(E_2)=2$. Since $E_1, E_2 \subset \ol {S_\mu x} $, $F$ is not
$S_\mu$-invariant but $\mu$-harmonic. This proves that (2) $\Ra$
(3).

Let $X$ and $Y$ be closed $S_\mu$-invariant subsets of $\ol {S_\mu
x}$ for $x\in E$.  Then $X$ and $Y$ contain $S_\mu$-minimal subsets.
This proves that (3) $\Ra $ (4). \eo

\end{subsection}

\end{section}

\begin{section}{Examples}

\bex\label{ex1} Consider the following linear action of $\R ^2$ on
$\R ^4$ given by
$$(t,s) \mapsto \pmatrix {1&0&0&0\cr 0&1&0&0\cr 0&0&e^{t-s}&0\cr
0&0&0&e^{s-t}}$$ for $s,t\in \R$.  Then the orbit of $v=(1,1,1,0)$
is $\{ (1,1,e^{t-s},0) \mid t,s \in \R \}$.  Let $\mP ^3(\R )$ be
the corresponding projective space and $\pi \colon \R ^4 \setminus
\{ 0 \} \to \mP ^3(\R )$ be the canonical projection.  We now claim
that the closure of the orbit of $\pi (v)$ has two disjoint
invariant sets.

Choose sequences $(t_n)$ and $(s_n)$ in $\R$ such that $t_n-s_n \to
-\infty$.  Then orbit closure of $v$ contains $(1,1,0,0)$.

Choose sequences $(t_n)$ and $(s_n)$ in $\R$ such that $t_n-s_n\to
\infty$.  Then $e^{s_n-t_n}(1,1,e^{t_n-s_n},0) \to (0,0,1,0)$.  This
implies that $\pi (0,0,1,0)$ is in the closure of the orbit of $\pi
(v)$.

It can easily be seen that $(1,1,0,0)$ is an invariant vector in $\R
^4$ and $\pi (0,0,1,0)$ is an invariant point in $\mP ^3(\R )$.
Thus, the closure of the orbit of $\pi (v)$ has two disjoint
invariant sets $\pi (1,1,0,0)$ and $\pi (0,0,1,0)$.

But since $\R ^2$ is abelian, $\R ^2$ has Liouville property and
hence any action of $\R ^2$ also has Liouville property.  \eex

The following example produces a measure on the $ax+b$-group that
does not have Choquet-Deny but has Liouville for its action on $\R
^2$.

\bex Let $G= \{ \pmatrix {t^2 & a \cr 0& t} \mid t>0,~~ a\in \R \}$.
Then $G$ is a solvable group and $G$ is the $ax+b$-group.  Any
measure $\mu$ on $G$ supported on $\{ \pmatrix {t^2 & a \cr 0& t }
\mid 0<t<1/5,  ~~ |a|<1/5 \}$ satisfies the condition in Theorem
\ref{gc}.  Hence $\mu$ has Liouville on $\R ^2$.  But $\mu$ itself
is not Liouville, that is, there are non-constant continuous bounded
$\mu$-harmonic functions on $G$ - this could be seen from section
5.1.2 of \cite{Ba}. \eex

The next example provides measures $\mu$ on $GL(V)$ that has
Liouville property on $V$ but neither $\mu$ nor $\check \mu$ has
Choquet-Deny.

\bex Assume that $V$ has dimension at least two.  Let $0<a <1$ and
$\mu$ be a probability measure on $GL(V)$ supported on $\{ g\in
GL(V) \mid ||g|| \leq a \}$. Then by Theorem \ref{gc}, $\mu$ has
Liouville on $V$. If $G_\mu$ is nonamenable , then neither $\mu$ or
$\check \mu$ can have Liouville property (or Choquet-Deny).  It may
be noted that $G_\mu$ is nonamenable if the support of $\mu$ is $\{
g\in GL(V) \mid ||g|| \leq a \}$.  \eex

\end{section}

\bigskip\medskip
\advance\baselineskip by 2pt
\begin{tabular}{ll}
C.\ R.\ E.\ Raja \\
Stat-Math Unit \\
Indian Statistical Institute (ISI) \\
8th Mile Mysore Road \\
Bangalore 560 059, India.\\
creraja@isibang.ac.in
\end{tabular}


\begin{thebibliography}{Abc}

\footnotesize


\bibitem {Ba} M. Babillot, An introduction to Poisson boundaries of Lie
groups.  Probability measures on groups: recent directions and trends,
1--90, Tata Inst. Fund. Res., Mumbai, 2006.

\bibitem {BQ} Y. Benoist and J-F. Quint, Random walks on projective
spaces, Preprint.

\bibitem {Da} S. G. Dani, On ergodic quasi-invariantmeasures of group
automorphism, Israel J. Math. 43 (1982), 62–74.

\bibitem {FeR} R. Feres and E. Ronshausen, Harmonic functions over group actions.
Geometry, rigidity, and group actions, 59–71, Chicago Lectures in
Math., Univ. Chicago Press, Chicago, IL, 2011.

\bibitem {Fu1} H. Furstenberg, Boundary theory and stochastic processes
on homogeneous spaces.  Harmonic analysis on homogeneous spaces (Proc.
Sympos. Pure Math., Vol. XXVI, Williams Coll., Williamstown, Mass., 1972),
193--229. Amer. Math. Soc., Providence, R.I., 1973.

\bibitem {Fu2} H. Furstenberg, A note on Borel's density
theorem, Proc. Am. Math. Soc. 55 (1976), 209-212.

\bibitem {GR} Y. Guivarc'h and C. R. E. Raja, Recurrence and ergodicity of
random walks on linear groups and on homogeneous spaces, Ergodic
Theory and Dynamical Systems (2012), 32, 1313-1349.

\bibitem {He} H. Heyer, Probability measures on locally compact groups.
Ergebnisse der Mathematik und ihrer Grenzgebiete, Band 94. Springer-Verlag,
Berlin-New York, 1977.

\bibitem {Hu} J. E. Humphreys, Introduction to Lie algebras and
representation theory. Graduate Texts in Mathematics, Vol. 9. Springer-Verlag,
New York-Berlin, 1972.

\bibitem {JaR} W. Jaworski and C. R. E. Raja, The Choquet-Deny
theorem and distal properties of totally disconnected locally
compact groups of polynomial growth, New York J. Math. 13 (2007),
159-174.

\bibitem {KaV} V. A. Kaimanovich and A. M. Vershik, Random walks on discrete
groups: boundary and entropy, Ann. Probab. 11 (1983), 457–490

\bibitem {Ro} J. Rosenblatt, Ergodic and mixing random walks on locally
compact groups, Math. Ann. 257 (1981), 31–42.

\bibitem {Ro2} M. Rosenblatt, Limits of convolution sequences of
measures on a compact: topological semigroup, J. Math. Mech. (1960),
293—306.

\bibitem {Wi} G. Willis, Probability measures on groups and some related
ideals in group algebras, J. Funct. Anal. 92 (1990), 202–263.

\end{thebibliography}
\end{document}